\documentstyle[12pt]{amsart}
\begin{document}
\title{Biharmonic submanifolds in a Riemannian manifold with non-positive curvature}
\author{Nobumitsu Nakauchi and Hajime URAKAWA${}^{\ast}$ }
\title[Biharmonic submanifolds]
{Biharmonic submanifolds in a Riemannian manifold with non-positive curvature}

\thanks{$^{\ast}$\,\,Supported by the Grant-in-Aid for the Scientific Research, 
(C), No. 21540207, 
Japan Society for the Promotion of Science.}
\dedicatory{}
\maketitle
\begin{abstract} 
In this paper, we show that, for every biharmonic 
submanifold $(M,g)$ of a Riemannian manifold $(N,h)$
with non-positive sectional curvature, if 
$\int_M\vert \eta \vert^2 v_g<\infty$, 
then $(M,g)$ is minimal in $(N,h)$, i.e., $\eta\equiv0$,
where $\eta$ is the mean curvature tensor field 
of $(M,g)$ in $(N,h)$. 
This result gives an affirmative answer 
under the condition 
$\int_M\vert \eta \vert^2 v_g<\infty$
to 
the following {\em generalized B.Y. Chen's conjecture}: 
every biharmonic submanifold of a Riemannian manifold 
with non-positive sectional curvature must be minimal. 
The conjecture turned out false 
in case of an incomplete Riemannian manifold
$(M,g)$
by  
a counter example of Y-L. Ou and L. Tang 
\cite{OT}
\end{abstract}
\numberwithin{equation}{section}
\theoremstyle{plain}
\newtheorem{df}{Definition}[section]
\newtheorem{th}[df]{Theorem}
\newtheorem{prop}[df]{Proposition}
\newtheorem{lem}[df]{Lemma}
\newtheorem{cor}[df]{Corollary}
\newtheorem{rem}[df]{Remark}
\section{Introduction and statement of results.} 
This paper is an extension of our previous paper \cite{NU} 
to biharmonic submanifolds of any co-dimension of 
a Riemannian manifold of pon-positive curvature. 
Let us consider an isometric immersion 
$\varphi:\,(M,g)\rightarrow (N,h)$
of a Riemannian manifold $(M,g)$ of dimension $m$ 
into another Riemannian manifold 
$(N,h)$ of dimension $n=m+p$ $(p\geq 1)$. We have
$$
\nabla^N_{\varphi_{\ast}X}\varphi_{\ast}Y=\varphi_{\ast}(\nabla_XY)+
B(X,Y),
$$
for vector fields $X$ and $Y$ on $M$, 
where $\nabla$, $\nabla^N$  
are the Levi-Civita connections of 
$(M,g)$ and $(N,h)$, and $B:\,\Gamma(TM)\times\Gamma(TM)\rightarrow \Gamma(TM)^{\perp}$ is the second fundamental form 
of the immersion $\varphi$ 
corresponding to the decomposition: 
$$
T_{\varphi(x)}N=d\varphi(T_xM)\oplus d\varphi(T_xM)^{\perp}\quad (x\in M), 
$$
respectively. 
Let 
$\eta$ be the mean curvature vector field along $\varphi$ 
defined by 
$
\eta=\frac{1}{m}\sum_{i=1}^mB(e_i,e_i), 
$
where $\{e_i\}_{i=1}^m$ is a local orthonormal frame on $(M,g)$. 
Then, the {\bf generalized B.Y. Chen's conjecture} (cf. \cite{CMO}, \cite{CMP}, \cite{C}, \cite{C2}, 
\cite{On}, \cite{Ou1}, \cite{OT}) 
is that: 
\vskip0.3cm\par
{\it For an isometric immersion
$\varphi:\,(M,g)\rightarrow (N,h)$, assume that 
the sectional curvature of $(N,h)$ is non-positive. 
If $\varphi$ is biharmonic (cf. See Sect. 2), then, it is minimal, i.e., $\eta\equiv0$. }
\vskip0.6cm\par
In this paper, we will show
\begin{th}
Assume that $(M,g)$ is a complete Riemannian manifold 
of dimension $m$ 
and $(N,h)$ is a Riemannian manifold of dimension 
$m+p$ $(p\geq 1)$ 
whose sectional curvature 
is non-positive. 
If $\varphi:\,  (M,g)\rightarrow (N,h)$ is 
 biharmonic
and satisfies that $\int_M\vert\eta\vert^2\,v_g<\infty$, then, 
$\varphi$ is minimal. 
\end{th}
\vskip0.6cm\par
In our previous paper \cite{NU}, we showed 
\begin{th}
Assume that $(M,g)$ is complete 
and 
the Ricci tensor ${\rm Ric}^N$ of $(N,h)$ satisfies that 
\begin{equation}
{\rm Ric}^N(\xi,\xi)\leq \vert A\vert^2.
\end{equation}
If $\varphi:\,(M,g)\rightarrow (N,h)$ is biharmonic and satisfies that 
\begin{equation}
\int_MH^2\,v_g<\infty,
\end{equation}
then, $\varphi$ has constant mean curvature,  i.e., $H$ is constant. 
\end{th} 
Notice that, in Theorem 1.2 
in case of codimension one, 
we only need the weaker assumption, 
non-positivity of the Ricci curvature of $(N,h)$ 
(\cite{OT}). 
On the other hand, in Theorem 1.1, 
we should 
treat with a complete submanifold of 
an arbitrary
co-dimension $p\geq 1$, 
and we need the stronger assumption 
non-positivity of the sectional curvature of $(N,h)$. 
In proving Theorem 1.1, 
the method of the proof of Theorem 1.2 
(\cite{NU}) does not 
work anymore. 
We should turn our mind, and have a different and 
very simple proof. 
Finally, our Theorem 1.1 
implies that 
{\em the generalized B.Y. Chen's conjecture holds true
under the assumption that $\int_M\vert \eta \vert^2\,v_g$ is finite and $(M,g)$ is complete}. 
\section{Preliminaries.}
\subsection{Harmonic maps and biharmonic maps}
In this subsection, we prepare general 
materials about 
harmonic maps and biharmonic maps of a 
complete Riemannian manifold into another Riemannian manifold (cf. \cite{EL}). 
\par
Let $(M,g)$ be an $m$-dimensional complete Riemannian manifold, and 
the target space $(N,h)$ is 
an $n$-dimensional Riemannian manifold. 
For every $C^{\infty}$ map $\varphi$ of 
$M$ into $N$. 
Let  $\Gamma(\varphi^{-1}TN)$ be 
the space 
of $C^{\infty}$ sections of the induced bundle 
$\varphi^{-1}TN$ of the tangent bundle $TN$ by $\varphi$.  
The {\it tension field} $\tau(\varphi)$ is 
defined globally on $M$ by 
\begin{equation}
\tau(\varphi)=\sum_{i=1}^mB(\varphi)(e_i,e_i)\in \Gamma(\varphi^{-1}TN),
\end{equation}
where the second fundamental form 
$B(\varphi)$ is defined by 
$$B(\varphi)(X,Y)=\nabla^N_{\varphi_{\ast}(X)}\varphi_{\ast}(Y)
-\varphi_{\ast}(\nabla_XY)$$
for $X,Y\in {\frak X}(M)$. 
Then, a $C^{\infty}$ map 
$\varphi:(M,g)\rightarrow (N,h)$ is {\it harmonic} 
if $\tau(\varphi)=0$. 
The {\it bitension field} $\tau_2(\varphi)$ is 
defined globally on $M$ by 
\begin{equation}
\tau_2(\varphi)=J(\tau(\varphi))
=\overline{\Delta}\tau(\varphi)-{\mathcal R}(\tau(\varphi)),
\end{equation}
where 
\begin{align*}
J(V)&:=\overline{\Delta}V-{\mathcal R}(V),\\
\overline{\Delta}V&
:={\overline{\nabla}}^{\ast}\,{\overline{\nabla}}V
=
-\sum_{i=1}^m
\{
{\overline{\nabla}}_{e_i}({\overline{\nabla}}_{e_i}V)
-{\overline{\nabla}}_{\nabla_{e_i}e_i}V
\},\\
{\mathcal R}(V)&:=\sum_{i=1}^m
R^N(V,\varphi_{\ast}(e_i))\varphi_{\ast}(e_i).
\end{align*}
Here, ${\overline{\nabla}}$ is the induced connection 
on the induced bundle $\varphi^{-1}TN$, and 
$R^N$ is the curvature tensor of $(N,h)$ 
(cf. \cite{H}) given by 
$$R^N(U,V)W=
[\nabla^N_{\,\,U},\nabla^N_{\,\,V}]W-\nabla^N_{\,\,[U,V]}W 
\quad (U,V,W\in {\frak X}(N)).$$  
\par
A $C^{\infty}$ map
$\varphi:(M,g)\rightarrow (N,h)$ is called to be 
{\it biharmonic} (\cite{CMP}, \cite{EL}, \cite{J}) if 
\begin{equation}
\tau_2(\varphi)=0.
\end{equation}
\vskip0.3cm\par
\subsection{Setting of isometric immersions}
In this sebsection, we prepare fundamental materials of general facts on 
isometric immersions (cf. \cite{KN}). 
Let $\varphi$
be an isometric immersion of an 
$m$-dimensional Riemannian into 
an $(m+p)$-dimensional Riemannian manifold $(N,h)$. 
Then, 
the induced bundle 
$\varphi^{-1}TN$ of the tangent bundle $TN$ of $N$ by 
$\varphi$ 
is decomposed into 
the direct sum: 
\begin{equation}
\varphi^{-1}TN=\tau M\oplus \nu M,
\end{equation}
where
$\varphi^{-1}TN=\cup_{x\in M}T_{\varphi(x)}N$, 
$\tau M=d\varphi(TM)=\cup_{x\in M}d\varphi(T_xM)$, and 
$\nu M=\cup_{x\in M}d\varphi(T_xM)^{\perp}$ is the normal bundle. 
 For the induced connection 
 $\overline{\nabla}$ on $\varphi^{-1}TN$ of the Levi-Civita 
 connection $\nabla^N$ of $(N,h)$ by $\varphi$, 
 $\overline{\nabla}_X(d\varphi(Y))$ is decomposed  corresponding to $(2.4)$ 
 as 
 \begin{equation}
 \overline{\nabla}_X(d\varphi(Y))=d\varphi(\nabla_XY)+B(X,Y)
 \end{equation}
 for all $C^{\infty}$ vector fields $X$ and $Y$ on $M$. 
 Here, $\nabla$ is the Levi-Civita connection of $(M,g)$ and 
 $B(X,Y)$ is the second fundamental form 
 of the immersion $\varphi:(M,g)\rightarrow (N,h)$. 
 \par
 Let $\{\xi_1,\cdots,\xi_p\}$ be a local unit normal vector fields 
 along $\varphi$ that are orthogonal at each point, and 
 let us decompose $B(X,Y)$ as 
 \begin{equation}
 B(X,Y)=\sum_{i=1}^pb^i(X,Y)\,\xi_i, 
 \end{equation} 
 where 
 $b^i(X,Y)$ $(i=1,\cdots,p)$ are the $p$ second fundamental forms of $\varphi$. 
For every 
$\xi\in \Gamma(\nu M)$, 
$\overline{\nabla}_X\xi$, denoted also by 
$\nabla^N_X\xi$  is decomposed correspondingly to (2.4) 
into 
\begin{equation}
\nabla^N_X\xi=-A_{\xi}(X)+\nabla^{\perp}_X\xi,
\end{equation}
where $\nabla^{\perp}$ is called the normal connection of 
$\nu M$. 
The linear operator $A_{\xi}$ 
of $\Gamma(TM)$ into itself, called the shape operator with respect to 
$\xi$, satisfies that 
\begin{equation}
\langle A_{\xi}(X), Y\rangle =\langle B(X,Y),\xi\rangle
\end{equation}
for all $C^{\infty}$ vector fields $X$ and $Y$ on $M$. 
Here, we denote the Riemannian metrics $g$ and $h$ simply by 
$\langle\,\cdot,\,\cdot\rangle$.  
\vskip0.3cm\par
We denote the tension field $\tau(\varphi)$ of an isometric immersion 
$\varphi:\,(M,g)\rightarrow (N,h)$ 
as 
\begin{align}
\tau(\varphi)&=
{\rm Trace}_g(\widetilde{\nabla}d\varphi)
=\sum_{i=1}^mB(e_i,e_i)\nonumber\\
&=\sum_{k=1}^p({\rm Trace}_gb^k)\,\xi_k\nonumber\\
&=m\sum_{k=1}^pH_k\,\xi_k\nonumber\\
&=m\,\eta,
\end{align}
where 
$\widetilde{\nabla}$ is the induced connection 
on $TM\otimes \varphi^{-1}TN$, 
$H_k:=\frac{1}{m}{\rm Trace}_gb^k
=\frac{1}{m}
{\rm Trace}_g(A_{\xi_k})$ 
$(k=1,\cdots,p)$, and 
$\eta:=\sum_{k=1}^pH_k\,\xi_k$ is the mean curvature vector field of $\varphi$.  
Let us recall that 
$\varphi:\,(M,g)\rightarrow (N,h)$ is {\em minimal} 
if $\eta\equiv0$. 
\vskip0.6cm\par
\section{Proof of Main Theorem}
Assume that
$\varphi:\,(M,g)\rightarrow (N,h)$ is a biharmonic immersion. 
Then, 
since (2.9): $\tau(\varphi)=m\,\eta$, 
the biharmonic map equation
\begin{align}
\tau_2(\varphi)=\overline{\Delta}(\tau(\varphi))-{\mathcal R}(\tau(\varphi))=0
\end{align}
is equivalent to that 
\begin{equation}
\overline{\Delta}\eta
-\sum_{i=1}^mR^N(\eta,d\varphi(e_i))d\varphi(e_i)=0.
\end{equation}
Take any point $x_0$ in $M$, and 
for every $r>0$, let us consider 
the follwoing cut-off function 
$\lambda$ on $M$: 
$$
\left\{
\begin{aligned}
&0\leq\lambda(x)\leq 1\quad (x\in M),\\
&\lambda(x)=1\qquad\quad (x\in B_r(x_0)),\\
&\lambda(x)=0\qquad\quad (x\notin B_{2r}(x_0))\\
&\vert\nabla \lambda\vert\leq \frac{2}{r}\,\,\,\,\quad\quad (\text{on}\,\,\,M), 
\end{aligned}
\right.
$$ 
where $B_r(x_0):=\{x\in M:\,d(x,x_0)<r\}$ and 
$d$ is the distance of $(M,g)$. 
In both sides of $(3.2)$, taking 
inner product with 
$\lambda^2\,\eta$, and integrate them over $M$, 
we have 
\begin{equation}
\int_M
\langle \overline{\Delta}\eta,\lambda^2\,\eta\rangle\,v_g
=\int_M\sum_{i=1}^m
\langle R^N(\eta,d\varphi(e_i))d\varphi(e_i),\eta\rangle\,\lambda^2\,v_g.
\end{equation}
Since the sectional curvature of $(N,h)$ is non-positive, $h(R^N(u,v)v,u)\leq 0$ for all 
tangent vectors $u$ and $v$ at $T_yN$ $(y\in N)$, 
the right hand side of $(3.3)$ is non-positive, i.e., 
\begin{equation}
\int_M\langle\overline{\Delta}\eta,\lambda^2\,\eta\rangle\,v_g\leq 0.
\end{equation}
On the other hand, the right hand side coincides with 
\begin{align}
\int_M
\langle\overline{\nabla}\eta,
\overline{\nabla}(\lambda^2\,\eta)\rangle\,v_g
&=\int_M
\sum_{i=1}^m
\langle\overline{\nabla}_{e_i}\eta,
\overline{\nabla}_{e_i}(\lambda^2\,\eta)\rangle\,v_g\nonumber\\
&=
\int_M\lambda^2\,\sum_{i=1}^m\vert\overline{\nabla}_{e_i}\eta\vert^2\,v_g\nonumber\\
&\quad+2\int_M
\sum_{i=1}^m
\lambda\,(e_i\lambda)\,\langle\overline{\nabla}_{e_i}\eta,\eta\rangle\,v_g,
\end{align}
since 
$\overline{\nabla}_{e_i}(\lambda^2\,\eta)
=\lambda^2\overline{\nabla}_{e_i}\eta
+2\lambda(e_i\lambda)\,\eta$. 
Therefore, we have 
\begin{align}
\int_M
\lambda^2\,\sum_{i=1}^m\vert\overline{\nabla}_{e_i}\eta\vert^2\,v_g
\leq 
-2\int_M
\sum_{i=1}^m\langle\lambda\,\overline{\nabla}_{e_i}\eta,(e_i\lambda)\,\eta\rangle\,v_g. 
\end{align}
\par
Now apply with $V:=\lambda\,\overline{\nabla}_{e_i}\eta$, and $W:=(e_i\lambda)\,\eta$, 
to 
Young's inequality: for all 
$V$, $W\in \Gamma(\varphi^{-1}TN)$ and 
$\epsilon>0$, 
$$
\pm2\,\langle V,W\rangle\leq \epsilon\,\vert V\vert^2
+\frac{1}{\epsilon}\,\vert W\vert^2,
$$
the right hand side of (3.6) 
is smaller than or equal to 
\begin{equation}
\epsilon\,\int_M
\lambda^2\,\sum_{i=1}^m\vert\overline{\nabla}_{e_i}\eta\vert^2\,v_g
+\frac{1}{\epsilon}\,\int_M
\vert\eta\vert^2\,\sum_{i=1}^m\vert e_i\lambda\vert^2\,v_g.
\end{equation}
By taking $\epsilon =\frac12$, we obtain
\begin{equation*}
\int_M\lambda^2\,\sum_{i=1}^m
\vert\overline{\nabla}_{e_i}\eta\vert^2\,v_g
\leq \frac12\,\int_M\lambda^2\,\sum_{i=1}^m
\vert\overline{\nabla}_{e_i}\eta\vert^2\,v_g
+2\,\int_M\vert\eta\vert^2\,
\sum_{i=1}^m\vert e_i\lambda\vert^2\,v_g.
\end{equation*}
Thus, we have 
\begin{align}
\int_M
\lambda^2\,\sum_{i=1}^m
\vert\overline{\nabla}_{e_i}\eta\vert^2\,v_g
&\leq 4\,\int_M\vert\eta\vert^2\,\sum_{i=1}^m\vert e_i\lambda\vert^2\,v_g\nonumber\\
&\leq\frac{16}{r^2}\,\int_M\vert\eta\vert^2\,v_g<\infty.
\end{align}
Since $(M,g)$ is complete, we can tend $r$ to infinity, and then 
the left hand side goes to $\int_M\sum_{i=1}^m\vert\overline{\nabla}_{e_i}\eta\vert^2\,v_g$, 
we obtain 
\begin{equation}
\int_M\sum_{i=1}^m\vert\overline{\nabla}_{e_i}\eta\vert^2\,v_g\leq 0.
\end{equation}
Thus, we have $\overline{\nabla}_X\eta=0$ 
for all vector field $X$ on $M$. 
\par
Then, we can conclude that $\eta\equiv0$. 
For, applying  $(2.7)$: 
$$
\overline{\nabla}_X\xi_k=-A_{\xi_k}(X)+\nabla^{\perp}_X\xi_k,
$$
to $\eta=\sum_{k=1}^pH_k\xi_k$, 
we have 
\begin{equation}
0=\overline{\nabla}_X\eta
=
-A_{\eta}(X)
+\nabla^{\perp}_X\eta,
\end{equation}
which implies 
that, for all vector field $X$ on $M$, 
\begin{equation}
\left\{
\begin{aligned}
&A_{\eta}(X)=0,\\
&\nabla^{\perp}_X\eta=0.
\end{aligned}
\right.
\end{equation}
by comparing the tangential 
and normal components. 
Then, by the first equation of (3.11), we have 
\begin{align}
\langle B(X,Y),\eta\rangle
=\langle A_{\eta}(X),Y\rangle
=0,
\end{align}
for all vector fields $X$ and $Y$ on $M$. 
This implies that $\eta\equiv 0$ 
since $\eta=\frac{1}{m}\sum_{i=1}^mB(e_i,e_i)$.
\qed
\vskip2cm\par       

\vskip0.6cm\par
Graduate School of Science and Engineering, 
\par Yamaguchi University,
Yamaguchi, 753-8512, Japan.
\par
{\it E-mail address}: nakauchi@@yamaguchi-u.ac.jp
\vskip0.8cm\par
Division of Mathematics, Graduate School of Information Sciences, 
\par Tohoku University, 
Aoba 6-3-09, Sendai, 980-8579, Japan.
\vskip0.1cm\par
{\it Current Address}: 
Institute for International Education, 
\par Tohoku University,
Kawauchi 41, Sendai, 980-8576, Japan.
\par
{\it E-mail address}: urakawa@@math.is.tohoku.ac.jp
\end{document}